\documentstyle[12pt,psfig]{article}
\begin{document}

\setcounter{page}{1}

  \hspace{16.6em}
 {\Large   SUNY-NTG  98-02
  }
\vspace*{0.8cm}

  \begin{center}
{\bf  \large
     Valuation of path-dependent American options
             \\
    using a Monte Carlo approach
             \\
       }
\vspace*{0.4cm}\ \\
           H.\ Sorge
 \footnote{
  E-mail: sorge@alba.physics.sunysb.edu
          }
\vspace*{0.4cm}\ \\
  Department of Physics and Astronomy,
\vspace*{0.4cm}\ \\
  State University of New York at Stony Brook, NY 11794-3800

  \end{center}

\begin{abstract}
It is shown how to obtain accurate values for American options using Monte
Carlo simulation. The main feature of the novel  algorithm consists of
tracking the boundary between  exercise and hold regions via optimization of
a certain payoff function. We compare  estimates  from simulation for some
types of claims  with  results from binomial tree calculations and find very
good agreement. The novel method allows to calculate so far untractable
path-dependent option values.
\end{abstract}                           

\newpage


 Pricing of options is an important area of research in the finance 
 community. The field has been pioneered by
 Black and Scholes (1973). They
   made a major breakthrough by deriving
 a formula for the price of any 
      contingent claim 
  which depends  on a  non-dividend-paying stock.
  Using the assumption of no arbitrage, they were able to show
  that the price of a derivative security  can be expressed as
   the expected value of its discounted payoffs.
   The expectation is taken under the assumption of a
   {\em risk-neutral} evolution of the value of the underlying security.
  Merton  (1973)  generalized these ideas to situations in
   which the interest rates themselves fluctuate in time.
  Due to the complexity of the underlying dynamics, numerical
    methods have become increasingly popular in modern finance.
   They are used for a variety of purposes, for instance
   valuation of securities and stress testing of portfolios.
   Analytical solutions 
   for problems in finance have been found only for rather special     
    cases.  For instance, 
    in order to gain the
    solution for  prices of European  options
   Black and Scholes  had to assume that the  evolution of the
    underlying asset price
    (more precisely its natural logarithm) 
    follows a so-called Wiener  process with  time-independent 
     volatility. 
   European options  can be exercised only at expiration date.
   In contrast, 
   American-style options which can be exercised during the whole lifetime
     of the option can be valued only numerically.
   Furthermore, numerical methods have to be employed if the
    dimension of the problem increases (more than one state variable)
    or in case of more realistic approaches to the
    stochastic process like for some no-arbitrage models of interest rate
    evolution.
    
   The Monte Carlo approach lends itself very naturally to the
   evaluation of security prices and interest rates.
   Schematically, it consists of the following steps.
   First, sample paths of the  state variables 
   (asset prices and so on) 
    over the relevant time horizon  are  generated.
   The cash flows of the securities on each sample path
   are evaluated, as determined by the structure of the security in
   question. The discounted  cash flows are averaged over the
   sample paths.
    The  Monte Carlo method  is  very flexible, since it does not depend
    much on the specific nature of the underlying stochastic process.
    Its accuracy is also independent on the dimensionality of the problem
     which is its dominant advantage over more traditional numerical
      integration methods.
    It is outside the scope of this paper to discuss the various
   facets of recent research on the use of  Monte Carlo in the
      finance area.
    The reader may consult the  concise recent reviews of
   Boyle, Broadie, and Glassermann (1997) on these topics.
    In particular, they report on recent progress 
      in developing more  {\em efficient}  Monte Carlo algorithms.
    Standard  Monte Carlo methods converge notoriously slow
     (with 1/$\sqrt{N}$, $N$  being the  number of sampled paths).
     A central goal of recent research activity  has therefore been the
     refinement of variance reduction techniques  
      (antithetic and control variates), 
       importance sampling 
     and low-discrepancy random number sequences.
   If the development in  sciences like physics is any guide,
      Monte Carlo applications in  finance will  become even  more
     important. With 
    computational prowess further  increasing 
       the  Monte Carlo method  may even shed    
     its `brute force' image from the past
    which was founded on its rather slow convergence property. 
   
   Boyle  (1977) was the first  to propose the use of
   Monte Carlo methods for option pricing 
   in the literature.
   Since the exercise date of European contingent claims is fixed,
   the mechanics of a price evaluation 
   of European options employing the Monte Carlo 
    method is rather straightforward.
 Let us begin by 
  collecting  the   notation for the most important variables
  which we are going to  use
  throughout the rest of the paper. In general, we follow the notation
  in Hull's book (1993), except that the
   present time for which the option price is calculated
  is designated to be ``0'' (without loss of generality).
\vspace{0.2cm}  

\noindent
\begin{minipage}[t]{0.30cm}
  \hspace{0.2cm}
\end{minipage} \
\begin{minipage}[t]{1.7cm}
  $S_t$
\end{minipage} \
\begin{minipage}[t]{0.35cm} =\end{minipage} \ 
\begin{minipage}[t]{13cm}
  the value of the state  (asset price, \ldots) at time $t$
   \\  underlying the derivative security
   \\ (which may be a single  variable or a vector), 
\end{minipage}
 \\[0.8em]
\begin{minipage}[t]{0.30cm}
  \hspace{0.2cm}
\end{minipage} \
\begin{minipage}[t]{1.7cm}
  $ S_0 $   
\end{minipage} \
\begin{minipage}[t]{0.35cm} =\end{minipage} \ 
\begin{minipage}[t]{13cm}
  the initial value  of the state variables, \\
\end{minipage}
 \\[0.8em]
\begin{minipage}[t]{0.30cm}
  \hspace{0.2cm}
\end{minipage} \
\begin{minipage}[t]{1.7cm}
   X
\end{minipage} \
\begin{minipage}[t]{0.35cm} =\end{minipage} \ 
\begin{minipage}[t]{13cm}
  strike price of the option,
\end{minipage}
 \\[0.8em]
\begin{minipage}[t]{0.30cm}
  \hspace{0.2cm}
\end{minipage} \
\begin{minipage}[t]{1.7cm}
   p
\end{minipage} \
\begin{minipage}[t]{0.35cm} =\end{minipage} \ 
\begin{minipage}[t]{13cm}
  path on which the underlying state  evolves in time,
\end{minipage}
 \\[0.8em]
\begin{minipage}[t]{0.30cm}
  \hspace{0.2cm}
\end{minipage} \
\begin{minipage}[t]{1.7cm}
   $G _t(p)$
\end{minipage} \
\begin{minipage}[t]{0.35cm} =\end{minipage} \ 
\begin{minipage}[t]{13cm}
  the value of the option {\em if} exercised at time $t$\\
   ($\max (S_t-X,0)$ for a simple Call and  \\
   $\max (X-S_t,0)$ for a Put),
\end{minipage}
 \\[0.8em]
\begin{minipage}[t]{0.30cm}
  \hspace{0.2cm}
\end{minipage} \
\begin{minipage}[t]{1.7cm}
  $r_t$
\end{minipage} \
\begin{minipage}[t]{0.35cm} =\end{minipage} \ 
\begin{minipage}[t]{13cm}
  instantaneous risk-less interest rate at time $t$, 
\end{minipage}
 \\[0.8em]
\begin{minipage}[t]{0.30cm}
  \hspace{0.2cm}
\end{minipage} \
\begin{minipage}[t]{1.7cm}
  $  \sigma $
\end{minipage} \
\begin{minipage}[t]{0.35cm} =\end{minipage} \ 
\begin{minipage}[t]{13cm}
  volatility (standard deviation) of state variable \\
     or their components respectively,
\end{minipage}
 \\[0.8em]
\begin{minipage}[t]{0.30cm}
  \hspace{0.2cm}
\end{minipage} \
\begin{minipage}[t]{1.7cm}
  T
\end{minipage} \
\begin{minipage}[t]{0.35cm} =\end{minipage} \ 
\begin{minipage}[t]{13cm}
  the lifetime of the option contract (expiration date),
\end{minipage}
 \\[0.8em]
\begin{minipage}[t]{0.30cm}
  \hspace{0.2cm}
\end{minipage} \
\begin{minipage}[t]{1.7cm}
 $V_t$
\end{minipage} \
\begin{minipage}[t]{0.35cm} =\end{minipage} \ 
\begin{minipage}[t]{13cm}
  the option value at time $t$. 
\end{minipage}
     
 We may easily estimate the 
 European option value  from a Monte  Carlo 
  path sample  for the variables related to the derivative security.
 These  paths  need to be generated 
    with probabilities determined by  the
     underlying stochastic process.
    The option value  at  the expiration date  $T$
    has been  specified  in the  option contract. 
    In general, it may 
    depend on the path   of the relevant state variables 
    between closing the contract
    and  expiration date.  Each  pathwise 
    final option value from the simulation is  discounted
    backward in time to determine its value at time 0.
    Since we do not know at time $t$=0 yet on which path 
     the  underlying variables will evolve, 
     the average over all paths is taken to estimate the present option
     value: 
   \begin{equation}
     \label{eq_europt_mc}
     V_0  \approx N^{-1} \sum _{p=1}^N 
               \exp( -r T )  V_{T}(p)
   \end{equation}    
   The discount factor includes the risk-free interest rate $r$,
     possibly averaged over the  lifetime  $T$ of the option. 
   For simplicity, it is assumed in eq.~(\ref{eq_europt_mc}) that
     all sampled paths have equal probability. Generalization
   to a situation in which paths are characterized by different
   probabilities which may arise e.g.\ in the context of
     Monte Carlo  optimization by importance sampling is straightforward.

  The pricing of  options 
     which may be exercised prior to the expiration date 
   by  Monte Carlo  simulation  is  more involved.  
   Here, the owner holds the right to exercise the option
    at several (Bermudan) or 
   possibly  infinitely many (American) `decision dates'.
    Many types of American contingent claims trade on exchanges
    and in the over-the-counter market.
     Examples include options, swaptions,
     binary options and Asian options. 
   It has been shown  by Roll  (1977), Geske  (1979), and
   Whaley  (1981) that exercise of
   {\em call } options which give the
    owner the right to purchase some underlying asset at the agreed-on
     strike price  is usually unfavorable
     before  expiration date, except close to ex-dividend
    dates.
    The situation is much less clear-cut in case of  {\em put } options
    which give the owner the right to sell the asset at the strike price.
    In these situations the `naive' Monte Carlo is encountering 
      unsurmountable difficulties. 
    In order to decide whether it is more favorable at some intermediate
     time to hold  or exercise the option the owner needs to compare
    the expected payoff in the two cases. The maximum of the two
     forms the option value at that time which 
      -- after discounting and taking the expectation value --
     leads to the following expression for the option
     at closing date:
   \begin{equation}
     V_0  =     \left\langle   
                 \exp ( -r {\cal T} ) V_{ \cal T} 
                               \right\rangle
                \quad .
   \end{equation}
  The angular brackets denote the expectation value with respect
   to the (risk-free) probability measure. 
    The path-dependent times $ {\cal T} $  are the  so-called optimum
    stopping times which may be any of the decision dates. 
  
   A single path  provides  clearly
   insufficient information  to evaluate the  option  value 
   in case of non-exercise     
   at any of the decision dates.
   The insufficiency of the  naive Monte Carlo method to deal
    with the optimum stopping problem has lead 
     some authors like Hull (1993) to the claim in
   the literature that ``Monte Carlo  simulation can only be used
   for European-style  options".
    On the other side, 
   evaluation of
   American-style Options via Monte Carlo simulation   
   has found already some consideration in the literature
   as reviewed 
   Boyle, Broadie, and Glassermann (1997).
   Their common denominator is a `clever' estimate of the option prices
    on the decision dates by bundling  subsets of
   paths. These authors noted 
    that the suggested algorithms cannot be considered
  satisfactory yet. Either some  approximations may lead to
  uncontroled errors affecting the simulation results or the required large 
  computational effort limits their applicability, e.g.\ to a small
  number of  exercise dates.
        
   In this paper we are suggesting a completely different  strategy
    to calculate  American-style  put and call option prices via
    Monte Carlo  simulation. 
     During the sampling of the paths  we do {\em not} attempt
      to estimate  the expectation of the payoff 
      if the owner continues  to  hold   the  option.
       Instead, we use the sampled paths to
    evaluate the boundary between the early-exercise region     
    and the hold region   in the space of variables entering into the
    option contract. The location of this boundary
      is the crucial piece of information whose knowledge
     allows  the straightforward use of the  Monte Carlo procedure
      for  option price estimation.
     We treat the valuation of American options as an
      {\em optimization problem} of a certain payoff function which
       depends on the  set of sampled paths. 
       This function depends also on the
        exercise policy, i.e.\ the
        boundary between  the early-exercise region
      and the complementary hold region  (YES-NO boundary in the following).
     Maximation of the  payoff function provides an estimator of
   the YES-NO boundary which gets arbitrarily close 
      to the true location of the boundary
       with increasing  number of sampled paths. 
       For the plain vanilla put and call options, 
      the boundary is  rather simple at  any exercise date,  a  point
       for one state variable, a line for a  two-dimensional 
        space of state variables and so on.
     After the boundary has been estimated the 
      price  estimation proceeds as for European options, because
     the option prices  for points on the boundary are known.
      The only difference to price estimation of European options as in
     eq.~(\ref{eq_europt_mc}) is  that -- concerning paths crossing the
      YES-NO boundary  -- 
      each pathwise option value is discounted starting from
     the point at which the path crosses the boundary
      for the first time:
   \begin{equation}
       \label{eq_amopt_mc}    
     V_0  \approx N^{-1} \sum _{p=1}^N 
               \exp( -r {\cal T}(p) ) \:  G_{{\cal T}(p)}(p)
                \quad .
   \end{equation}    
    Time ${\cal T}(p)$ denotes the earliest time at which the $p$-th path
    crosses the YES-NO boundary (if at all). Otherwise  ${\cal T}(p)$
    equals the time until expiration of the option ($T$).
    $ G_{{\cal T}(p)}(p)$ is the agreed-on payoff 
   from exercising the option
     at  time  ${\cal T}(p)$ under the assumption that  path $p$ represents
    the realized evolution of the underlying  state variables.
      Here we have tacitly assumed that the starting point lies
     in the NO region. 
      Of course, the algorithm is able to handle the 
       other   possibility -- immediate exercise -- as well.
     In this case calculation of the option price is trivial, however.
      
  In a later section of this paper we will demonstrate that
   the  novel  Monte Carlo   algorithm  achieves an accuracy
  in the determination of American option values which compares
   well with  corresponding calculations for European options.
   We compare the results for a spectrum of standard American Puts
      with results from  binomial tree calculations. 
   This method  originally introduced by
     Cox, Ross, Rubinstein  (1979) is widely  used for option
    valuation.
   Here we assume an ideal Wiener-type  stochastic
      process for a single state variable, because a binomial tree
   is well suited to provide very accurate results in this situation.
  Of course, this comparison serves mostly illustrative purposes
   to show that the method works at finite $N$. 
  Lateron, we are going to present numerical results 
   for {\em path-dependent } American options. 
   We consider two cases,  Puts  on either the geometric or the 
    arithmetic average over the lifetime of the option.
   Assuming again an underlying ideal Wiener process 
    we may compare the  Monte Carlo estimations to results 
    from binomial tree calculations. While options on the geometric
    average can be priced by the tree method to arbitrary accuracy
     no such tree method is available for options contingent on the
      arithmetic average value.
    Some approximative method has been suggested 
    to express  the latter  option prices 
      by the corresponding prices in case of geometric averaging.
    Therefore we will compare the approximative formula
     to our more accurate procedure. 
     In general,
    it is envisioned that the novel  Monte Carlo  approach 
     and some variants will be
    applied for valuation problems   which 
   cannot or not easily be tackled by other methods.     
   Such problems encompass stochastic processes with Non-Markov 
   properties,stochastic volatilities,  multi-dimensional state variables
   and other path-dependent American options.

\newcounter{cnt_sect}
\setcounter{cnt_sect}{1}
\section*{ \Roman{cnt_sect} Valuation of 
        American-style options as an optimization problem}
      
 Let us assume in the following that the underlying state variables
  whose values determine the value of an option at exercise dates
  follow a  known stochastic process. 
  Generally, a stochastic process is specified by defining the state space,
  an index parameter (usually the time) and the dependence relation
   between the random variables. The latter is usually given in terms
 of a stochastic integro-differential equation containing some deterministic
  drift terms and  random changes in addition, notably diffusion or jump 
   processes. 
   Memory or so-called non-Markov effects
    may readily be included into  the transition probabilities.
  Non-Markov models have found recently considerable interest in the
 finance community, because the Heath, Jarrow and Morton  (1992) approach
  to interest rate evolution can be shown to possess 
  Non-Markovian features, in general.
    We assume  `smoothness' for the
   transition probability  density function
   which determines the
   probabilities -- based on a path's  evolution history -- 
   that the path  
   will comes close to any element of the state space at  some  time
     in the future. These
    are supposed to be continuous
     functions of the path variables. 
   A discrete set of sequentially 
   selected random  points in time is often called a
  ``chain". A path is the corresponding generalization 
    for a continuous time evolution. 
      
  At first, 
  we will restrict ourselves  to  options
   and stochastic processes 
  which result in the following simple property
   of the YES-NO boundary:
 \begin{itemize}
  \item[A1:]
   ``The {\em  boundary}
   between YES and NO regions   at any given decision time $t<T$
   is a simply connected hypersurface 
   in the state space with
   dimension of one unit less than
   the dimension of the  state space." 
 \end{itemize}
  
  As a consequence of the assumption,  there are no disconnected 
    pieces of    YES or  NO regions at a given time either.
   With the exception of some pathological cases, 
    assumption A1 will  be compatible only with
    option payoffs contingent on the values of
    the state variables at  exercise date. Instead of specifying
    these pathological  stochastic processes
   we simply require here:
 \begin{itemize}
  \item[A2:]
    ``The  payoff from exercising an option 
       is only a function of 
          the state variables at the exercise date
         but not at earlier times, i.e.\
        $G _t(p)= G_t (S_t,X)$."
 \end{itemize}
  
 In general, the topology of exercise  and hold  regions
  of  American  options
  may be more complicated than allowed by our
  assumptions. 
  Two  examples are
  some types of  barrier options and path-dependent options.
  Suppose the owner of a barrier option  receives only some
  payoff if a security stays within some specified range of values
  at exercise date.
  Therefore the option will not be exercised for `too small' and
   `too large' values of the underlying  which leads to
  two disconnected  NO regions.
 Another  perhaps more 
  interesting example of options not covered by the assumption
  stated above are American path-dependent options whose 
   valuation we will return to later in this section.
 Here the payoff  depends on the values which the underlying
 state variables have taken in the past, e.g.\ on the arithmetic mean
  within some time interval. Projected 
   into the space of state variables at any given time
  YES and NO region are completely scrambled. 
    It  may be optimal to exercise or hold  an option 
    for the {\em same} present 
    value of the state variables $S_t$ depending on the
   different histories of each path.
   Lateron, we will  show  that a  redefinition of the
   state space enables us to cover the case of
   path-dependent options by the same algorithm. 

  We should also mention that the assumption A2 restricts the
  memory   properties of the underlying stochastic process.
  They should not spoil the assumed
  property of the   YES-NO boundary in the state space. 
  This simply means that although the future evolution 
   of each path may depend
   on the past, the `average outcome' determining the
   decision whether to exercise the option at present time 
   is not influenced by the past but   of 
  Markovian nature.

   Assumption A1 allows us to define  coordinates  
  in  the  $d_S$ dimensional state space 
  such that
      any element $S_t$ of the state space  can be uniquely 
      characterized by some point of the  YES-NO boundary  
      $S_t^Y$
       and an additional coordinate $c_t$ 
 \begin{equation}
     \label{eq_point}
        S_t = (S_t^Y,c_t) 
 \end{equation} 
   in a  peculiar way.
    For any path passing through $S_t$ at present time
   the following inequalities hold  between the 
     payoff from exercising the option $ G_t $
     and the option value 
      $V_t^N  $ 
    if the option is being held at least until the next
    exercise date:
\begin{equation}
     \label{eq_coord}
   G_t (S_t,X)
     \left\{
           \begin{array}{ll}     
       > V_t^N   &
        \mbox{for $c_t < 0$}   \\
       < V_t^N  &
        \mbox{for $c_t> 0$}  \\
       = V_t^N   &
        \mbox{for  $c_t= 0$}
         \end{array}
         \right.        
 \end{equation} 
 The content of eq.~(\ref{eq_coord}) is to provide 
 convenient coordinates from which one can easily read off
  whether  point  $S_t$  is a boundary point ($c_t= 0$),
  in the 
 YES ($c_t\le 0$)  or in the NO ($c_t\ge 0$)  region.
  The knowledge of all $c_t$ or just their signs is tantamount
   to finding the YES-NO boundary itself. 
  Unfortunately,  
   a  straightforward application of the inequalities as expressed in
    eq.~(\ref{eq_coord}) to determine the boundary is prohibitive.
   Estimating the 
    expected payoff from holding the option 
     and comparing it to the payoff from exercising
    would require a `new simulation within the simulation'.
     Such a procedure becomes numerically awkward and would
     severely restrict the complexity of the problem to be
    tackled, concerning  e.g.\  the nature of the stochastic process 
     and the number of exercise dates. 
   
  We will proceed  now to estimate the location of the  YES-NO boundary
  based on a set of sampled paths. For this purpose
   it will be necessary to evolve the state variables backward in
  time. This feature of the novel Monte Carlo  approach which it shares with
   all  other numerical approaches to   option pricing
    like the binomial tree method  is caused by the
    boundary conditions of the problem. 
   The option price is a specified function of the state variables
   only on the expiration date. 
    Of course, the location of the  YES-NO boundary
     is also trivially determined at expiration time. 
    For a simple call or put option on a single asset 
     it is the point  $S_T=X$.
  
  The  continuous time process 
  will be presented as  a sequence of  time steps
   $\delta t_i , i=1,\ldots N_T$ which add up to the
    lifetime of the contract $\sum \delta t_i = T$. 
  The time steps are supposed to match the distance between decision times 
   exactly if a finite number of exercise opportunities has been specified
   in the option contract.
  For American options which give the owner the right of exercise
    at any time this is an approximation. However, choosing 
      `small'  time steps $ \delta t$ the results
     will be close to the continuum limit. 
  The current time $t$ for which the location of the
   YES-NO boundary  is being sought can be represented   by
   the integer index $i$ (with $t= i \cdot \delta t$
     for American options). 
  Since the location of the YES-NO boundary will be determined
   step by step moving backward in time, it is implicitly assumed
  from now on that the  YES-NO boundary has  already been determined 
  for later times $i$+1,$i$+2,\ldots ,$N_T$ (within the
  usual uncertainties due to the numerical imprecision).
  Therefore we are able to evaluate 
   for each path 
 \footnote{
  For the ease of reading
  we are  using calligraphic letters in this paper
  when dealing with functions of  paths  or path segments.
  }       
    its  payoff from optimal exercise at {\em later} times:
 \begin{equation} 
  \label{calG_i+1}
  {\cal G}_{i+1}(p_{>t})  := 
      \left\{
         \begin{array}{ll}
           G_t (S_{i+1}(p),X)  &  \mbox{for $S_{i+1}(p)$
                                   in  YES region}  \\
              &   \\
           \exp ( - \sum _{k=i+1}^{j-1} r_k \delta t_k ) 
                       \cdot
                      &       \\
                    \hspace{2.0cm}  G_t  (S_j(p),X)
                             \qquad
                         &  \mbox{for $S_{i+1}(p)$  
                                          in  NO region}
                \quad .
         \end{array}
      \right.
 \end{equation}
     The path-dependent payoff has been
     discounted to time $i+1$ in eq.~(\ref{calG_i+1}).
    $j$ represents the earliest time   (under the constraint $j>i$)  
  at which path $p$ crosses the boundary from the  NO region
   into the  YES region or equals $N_T$ if the path stays
   within the NO region for all times between $t$ and $T$.
   (For the latter case $ G_t  (S_j(p),X)$=0, of course.)
 We have used the subscript $>t$ for the argument $p$ of
 $ {\cal G} _i $ to emphasize that  the  payoff defined above 
    depends only on the values of the path variables for
    times  $> t$ but not on the past evolution of the state.
   So far, we are not able  to calculate the discounted pathwise
    payoff $ {\cal G}_{i}(p_{\ge t})$, i.e.\ for time $t$ . This would 
    require  the knowledge
    of the  YES-NO boundary at time $t$. 
However, 
   the  path-dependent payoffs  given in  eq.~(\ref{calG_i+1}) 
  enter into the expression for
   the value of the option at time $i$ 
    subject to the condition that
  the option is {\em not} immediately exercised (and, of course, 
     that the state variables have evolved along path $p$ so far):     
 \begin{eqnarray}
     \label{eq_opt_price_ih}
   V_t^N &=&
         \exp ( - r_i \delta t_i )
                  \int  D  p_{>t}  
                                \\  \nonumber
           &  &       \hspace{2.5cm}
               {\cal G}_{i+1}(p_{> t}) \:
                  \pi (p_{> t};p_{\le t},t)
                \quad .
 \end{eqnarray}   
  The integration  over all future path segments
   in  eq.~(\ref{eq_opt_price_ih})
  \begin{displaymath}
        \int D  p_{>t}  :=  \int  \prod _{l=i+1}^{N_T}
                               d ^{d_S} S_l
   \end{displaymath}
   is  an ordinary multiple integral (but a path integral
    in the continuous time limit).
     $ \pi (p_{> t};p_{\le t},t) $  
     denotes the  transition probability density function and
    describes the  probability (density) to evolve along
      the  path segment $p_{> t}$ in the future up to time   $T$ 
      depending on  the present 
    -- and previous (for  Non-Markov processes) --
       states of  path  $p_{\le t}$.

  Next we define a  continuous path $H$ in the state space
   by letting $c$ be a continuous parameter within 
    some range such that all
      $S'(c)=(S_t^Y,c) $ 
      are elements of the state space and vice versa.
  Let us assume for a moment that we  know already that 
       this line  intersects with the   YES-NO boundary
    but  we do not know which point it is. 
    (Somebody else has constructed the path $H$.)
      Can  we estimate 
      the value of the  YES-NO boundary point  using 
     the  information provided  by the path sample
     and obtain its correct value
       in the limit $N\rightarrow \infty$? 
    The answer is positive. We are able to construct 
     a function of the states  $S'(c) \in H$  whose only maximum
     is at  $S'(c)=(S_t^Y,0)$, i.e.\ the  boundary point, 
     for $N\rightarrow \infty$. 
   It is clear that our ability to define such a function 
     has important repercussions
   for the question how to {\em find } the  YES-NO boundary. 
    The  virtue of path $H$ is to provide us with 
     the knowledge that it has exactly  one 
     intersection point with the  YES-NO boundary.
    We do not necessarily need this information. 
   Instead, we may  search for the global maximum on
    continuous paths  connecting  points in the YES region
   with points in the NO region.
    According to statement A1  they will cut through the
    YES-NO boundary  at least once.
    (It will not hurt the search procedure, 
    if we accidentally hit the boundary several times
     reflected by several local maxima.)
   It is also rather straightforward to determine
    selected points which are in either the  YES or the NO
    region for sure. For instance, 
     points lie in the  NO  region if the payoff from exercising is
     zero.  Or the point with  zero value of the underlying 
        is certainly an element of the YES region for 
     a plain vanilla Put, 
      because the payoff cannot get any better.    
   
 In order to proceed with the proof of our conjecture we define
 a --  generally suboptimal  -- exercise policy
 for all paths  intersecting with $H$ at time $t$
  which depends on an arbitrary point $S'=(S_t^Y,c') \in H$. 
  We define a preliminary   ``YES'' region 
   (with the  quotation marks reminding us  that  the true
     YES boundary may be different) as to cover all paths
    with state variable $S_i(p)=(S_t^Y, c_i(p))$
     which fulfill  $ c_i(p) < c'$. 
   For instance, in case of a plain vanilla 
  put (call) option  on a single underlying asset we may declare that
  $S_i(p)$ lies in the  ``YES'' region if it fulfills  the condition
    $S_i(p) < $($>)$~$S'$.
   The pathwise payoff from the $S'$-dependent policy is
 \begin{equation}
     \label{eq_pisss}
       {\cal P} _i (p_{\ge t}, S') :=
         \left\{
           \begin{array}{ll}
               G_t (S_i(p),X)  &  \mbox{for $c_i(p) < c'$,}
                                     \\ 
           \exp ( - r_i \delta t_i ) 
                 \:   {\cal G}_{i+1}(p_{> t})  &  \mbox{else}
                \quad .
           \end{array}
         \right.
 \end{equation}   
  It is clear that the  pathwise payoffs will sometimes be larger 
  for `wrong' boundary points ($c'\ne 0$) than for the correct one.
  This reflects the  additional information 
  about the future which is encoded in the 
   path segment $p_{\ge t}$.
  The  essential point is now to sum over all  paths
   which removes this bias:
 \begin{equation}
     \label{eq_optim_mc}
     {\cal O} _t^N ( S') :=
       N(H)^{-1} \sum _{p}^{N(H)} {\cal P} _i (p_{\ge t}, S') 
        \quad  .
 \end{equation}   
  The function $ {\cal O} _t^N ( S')$  depends on time, 
    the sampled paths intersecting with $H$ 
  (of number $ N(H)$)   and
  on   $S'$. 
  Our central proposition 
   can be stated now as follows:      
  
 \begin{itemize}  
 \item[P1:]   
 ``Taking the limit  of infinitely many paths, 
  $ { O} _t ^\infty( S'):= 
    \lim _{N(H)\rightarrow \infty}  {\cal O} _t ^N$  
        has one extremum which is a  maximum
   within the  one-dimensional domain  $H$. Its 
    position is  the intersection point of $H$ with the
    YES-NO boundary   $S'$=$(S_t^Y,0)$."   
 \end{itemize}   

 In order to prove  proposition P1 we are going to rearrange the sum
 in  eq.~(\ref{eq_optim_mc}). 
  We  put all paths which have 
  essentially the same transition probability density function
   into the same (albeit infinitesimally small) `bin'.
   Here we need the assumption stated above that the
   transition probability density function  is a continuous
   function of the path variables.  
After `binning' and taking  the limit $N(H) \rightarrow \infty$
  eq.~(\ref{eq_optim_mc}) turns into 
 \begin{eqnarray}    
     \label{eq_optim_cont}
  { O} _t^\infty ( S') &= &   
       \int  D  p_{\le t} \; 
           \rho (p_{\le t})
                 \left( 
               \int D  p_{>t}  \: 
                   {\cal P} _i (p_{\ge t}, S')
                    \:   \pi (p_{> t};p_{\le t},t)
                                                \right) 
    \quad .
 \end{eqnarray}  
   The  transition probability density function
     $ \pi (p_{> t};p_{\le t},t) $   
       appears in eq.~(\ref{eq_optim_cont}), because 
       relative frequencies approach the  probabilities
       which are governing the random selection of paths according to
        Bernoulli's law of large numbers.
   The probability density
  $ \rho (p_{\le t})  $  
   to have a  path segment 
     close to $p_{\le t}$ needs  to
   be defined in accordance with the integration measure $ D  p_{ \le t}$
    and is normalized to one.

   The  quantity in brackets  $()$ 
   on the right hand side of  eq.~(\ref{eq_optim_cont}) 
 \begin{eqnarray}   
     D_t(S,S') &:= &
        \int D  p_{>t}  \:
                   {\cal P} _i (p_{\ge t}, S')
                    \:   \pi (p_{> t};p_{\le t},t)    
 \end{eqnarray}  
 can be readily evaluated using eqs.~(\ref{eq_opt_price_ih}) 
  and (\ref{eq_pisss}) to be
 \begin{eqnarray}   
     D_t(S,S')   &=&  
         \left\{
           \begin{array}{ll}
               G_t (S,X)  &  \mbox{for $c < c'$,}
                                     \\
               V_t^N (S)    &  \mbox{else}
                \quad .
           \end{array}
         \right.          
                      \\   \nonumber
 \end{eqnarray}  
 $D_t(S,S')$ coincides with the  option price  $V_t (S)$,
  except for  $c$ values in a certain range.
  Furthermore, on virtue of eq.~(\ref{eq_coord}) 
     $D_t(S,S')$ is always {\em smaller } 
       than  $V_t (S)$ for such chosen 
    $c$ values:
 \begin{eqnarray}
     \label{eq_Di}
   V_t (S) &>&  D_t(S,S')
         \left\{
           \begin{array}{ll}
            =  G_t (S,X)  &
                \mbox{for  $c \in (0,c')$ with $c'>0$,}
                                     \\
             =  V_t^N (S)  &
                \mbox{ $c \in (c',0)$ with $c'<0$}
           \end{array}
         \right.
         \\  \nonumber
     &=&  D_t(S,S')   \hspace{3em}   \mbox{ else }
                \quad .
 \end{eqnarray}   
   The conditions of eq.~(\ref{eq_Di})
    together with the positiveness of  the weight factors 
   $ \rho $  in eq.~(\ref{eq_optim_cont}) 
    guarantee  that 
 \begin{eqnarray}       
    { O}_t^\infty (c_a) & < & { O}_t^\infty (c_b)
                      \qquad \mbox{for $c_a < c_b \le 0$}  \\
                   \nonumber
    { O}_t^\infty (c_a) & < & { O}_t^\infty (c_b)
                      \qquad \mbox{for $c_a > c_b \ge 0$}  
      \quad .  
 \end{eqnarray}       
  This completes the
   proof of proposition P1.
 
  As stated above we may turn things around and construct the
     YES-NO boundary  from a set of lines 
   in state space which connect elements of the 
   YES and NO region.  The global 
     maximum of  function  ${\cal O} _t^N ( S') $  along each line
   is the supposed intersection point with the YES-NO boundary. 
   This process may be iterated until the desired precision
     is reached. 
  This standard optimization procedure becomes particularly simple
    in case of plain vanilla options. 
    Since the  YES-NO boundary is just a point at any exercise date,
    its location is already completely specified with one iteration step.
  The determination of the  YES-NO boundary  at time $t$ (index $i$)
    enables us to calculate the pathwise payoff at this time.
    Going  step-by-step  backward in time we may finally calculate
       the pathwise payoffs at initial time
    whose average determines the option value at contract time.
        
  The presented optimization technique may
   be used in the problem of valuation of 
  American {\em path-dependent} options as well.        
   Here the payoff from exercising the option depends on
   some function of the state variables in the past
   which we will denote as $\bar{S_t}(p_<)$
    in the following. 
   For instance,  $\bar{S_t}(p_<)$  may be the -- arithmetic or
    geometric -- average of any  state variable over some
    time window or the value at some specified time(s) of the past
    (look-back options).
   For simplicity, we restrict ourselves here to stochastic
    processes of Markovian nature whose parameters have been specified. 
   In this case the exercise decision  at any  time 
    before expiry is a function
   of just two  -- each possibly vector-valued -- 
    variables, the  values of the state variables
    and of the path-dependent variable  $\bar{S_t}(p_<)$, 
      both taken at the decision time.
   The difference between 
    standard  path-dependent and plain vanilla options 
    is just that the role of the current value of the
     underlying variable in  plain-vanilla  contracts
      is played by  $\bar{S_t}(p_<)$  in path-dependent 
     contracts, e.g.\ Max$(X-\bar{S_t}(p_<),0)$ for a Put.
   The topology of the exercise region for this kind of
    path-dependent options is  very simple but only in the
    space of state variables {\em enlarged by the additional
    dimension(s) of  $\bar{S_t}(p_<)$ }. 
    The  YES-NO boundary  in the $S_t$-$\bar{S_t}(p_<)$  space
      is again simply connected. Only its projection into the
    $S_t$ space leads to a scrambling of YES  and NO regions.
   This suggests immediately how the suggested optimization
    technique may be applied for the case of these  path-dependent
   options. We equip the  space of states with additional
   dimensions by declaring   $\bar{S_t}(p_<)$  to be one of the
   stochastic variables.
   Of course, the stochastic character of this  variable
    is rather funny. Its change with time is completely
    deterministic and of Non-Markovian type. However, these
   features are in accordance with the basic assumptions
     entering into the proof of  proposition P1. 
    The complications of such path-dependent option contracts
    just result in a higher dimension of the  space in which we
    have to track the  YES-NO boundary. 
        
   We are now in a position to compare our novel strategy
    for an extraction of  American-style  option prices
       from  Monte Carlo simulation to previously suggested ones.
    Boyle, Broadie, and Glassermann  discuss mainly three strategies
     in their review paper 
    to overcome the inherent limitations 
     of the naive  Monte Carlo procedure: the bundling algorithm
     invented by Tilley   (1993),
     the stratified state aggregation algorithm suggested by
     Barraquand and Martineau  (1995) 
      and  Broadie and  Glassermann's  (1995)
    algorithm based on simulated trees. 
     In Tilley's approach 
      bundles of paths which are `close' in state space at given time
      are considered to emanate from one single
     point for the purpose of option valuation.
      This introduces some error whose magnitude
     cannot be easily controled.
     The bias of the exercise decision will be very strong
      if a bundle  contains only   few paths. Such a situation is likely
     to emerge  in situations  
      in which there are many  exercise dates or several  state variables.
      Barraquand and Martineau  introduce another type of `bundling',
       however, in the payoff space.  
   Our approach shares similarities with  Tilley's  algorithm in the
   sense that 
    it  contains an implicit `bundling' of paths
    (called `binning') in the proof of the proposition P1.
     However, we have actually never to carry out this
     bundling in a numerical calculation, because we stick
     to an optimization of the  `inclusive' function
      $ { O} _t^\infty ( S')$ which sums up all bins.
     Broadie  and Glassermann proposed  an algorithm based on a
      `bushy tree' structure in state space
      which avoids partitioning of the state or payoff space.
      In their model
      many  branches emanate from each node.
       This process is replicated as often as  there are      
     decision times $d$ in the problem.
       The replication process effectively limits
         $d$ to small numbers (on the order of 4) in practical
      computations.  
      They also discuss the problem that the suggested  estimators
       of the option prices
       tend to be biased in `up' or `down' direction.            
        and converge only to the correct value 
     in the limit $N\rightarrow \infty$.
        This problem is of great relevance for  practical  
       calculations (which are always at finite $N$). They use
      the idea to calculate two estimators -- one biased low and
       the other high -- to obtain confidence limits from the
      simulations. Lateron, we are going to apply this idea
       in the framework of our approach.
      
\setcounter{cnt_sect}{2}
\section*{  \Roman{cnt_sect}  Simulation results for  American  options}
      
  In this paper we will only consider stochastic processes of
   rather simple nature, so-called ideal Wiener 
   or diffusion processes
   (for the $\log $ of any state variables).
   In physics these processes are usually called Brownian motion.
  Such  processes seem to be  most commonly
   used in valuation problems. Furthermore,  numerical
   techniques have been developed to calculate  prices for 
    plain vanilla   and other simple types of American options.
   In particular, the binomial tree method may provide very accurate
     values  and will be used as a benchmark here
    against which the results from the Monte Carlo simulations are
    being tested.
   In this section we will restrict ourselves to {\em Puts }
    on the current  value and on the -- geometric or arithmetic  -- 
     average  of  one single underlying security.
  An  ideal Wiener process  can be 
     characterized by the differential transition law 
 \begin{equation}
  \label{eq_Wiener}
   d\hat{S_t} = \mu  \cdot dt +  \sigma \sqrt{ dt } \cdot z  \quad ,
 \end{equation}   
    governing the change 
   of the state variable $\hat{S_t}:=\ln S_t$ within a  time interval
   $dt >0$.
  $z$ is a random number drawn from
 a normal distribution with mean 0 and variance 1.       
 $\mu$ is in general the drift velocity which may be replaced by
   the risk-free instantaneous interest rate $r$ using arbitrage
 arguments.
  Sometimes one needs to follow an  ideal Wiener process 
  only {\em backward } in time, e.g.\ in case of Monte Carlo estimation of
 plain vanilla options. Using the {\em bridge } construction
 --   start  $\hat{S_0}$ and final value $\hat{S_T}$  are
    known   --   
  a random value at intermediate time  may be chosen according to
 \begin{equation}
  \label{eq_Wiener_bw}
   d\hat{S_t} = 
            \sigma \sqrt{\frac{t}{t+dt}dt }\cdot z  
          - \frac{dt}{t+dt} (\hat{S}_{t+dt} - \ln S_0 )
           \quad .
 \end{equation}   
  After having fixed final and start values,  iterative use of 
    eq.~(\ref{eq_Wiener_bw}) may be employed to
   generate a time reverted copy    of the standard diffusive Wiener-type   
  motion.
  For the simulation results presented below we have 
   used a simple importance sampling scheme. 
  The exercise value  of strongly  out-of-the money options at expiry
  will be zero on  most paths. Therefore 
   the weight of the rare paths resulting in nonzero option value
    has been artificially enhanced in the path generation process.
  Of course, the enhancement factors are taken out again  in the calculation
  of the path-averaged prices.
      
  In its most simple form the   Monte Carlo algorithm  
    which is suitable for  valuation of  American-style  options
     consists of the following steps:
 \begin{description}
   \item[(1)]
    Generate path sample.
    \item[(2)] 
    Go step-by-step  backward in time starting at expiry.
   \item[(3)] 
     Track  the  YES-NO boundary 
       by employing the optimization of 
       function  ${\cal O} _t^N ( S') $
       along arbitrary paths connecting  YES and NO region
        (cf.\ eq.~(\ref{eq_optim_mc})).
   \item[(4)] 
     Evaluate  the pathwise payoffs at the earlier time
      analogous to eq.~(\ref{calG_i+1}). 
   \item[(5)] 
    Repeat the process until the contract time is reached. 
   \item[(6)]    
     Average  the pathwise  payoffs discounted to initial time
      over the whole  path sample
      (cf.\ eq.~(\ref{eq_amopt_mc})). 
 \end{description} 
  
  Let us discuss now some aspects related to the {\em finite }
   size of the path sample in practical calculations. 
 It is clear that step {\bf (6) } does not provide
  an unbiased estimator of the option price at any finite $N$
   but is biased upward. 
   The reason is that an  error is made  extracting  the
    YES-NO boundary from the finite size path sample.
    A  deviation of the extracted from the true boundary has its
   root in the  information about the mismatch between
  frequency distribution of the $N$ paths and underlying
    probabilities which enters into the optimization process.
   The  estimator becomes only asymptotically unbiased.
  Since we are not able to construct an unbiased estimator, we may
  get an estimate of the error introduced by the bias in step
     {\bf (6) } by constructing an estimator which is downward
   biased at finite $N$ but also asymptotically unbiased.
  It is very simple to construct such an estimator: 
 \begin{description}
   \item[(6')] 
   Generate a
   path sample independent from the first one and evaluate the
   American option prices using the previously calculated 
    YES-NO boundary.
 \end{description} 
  The reason that the prices  calculated from this path sample
   are biased downward mirrors the one for the upward bias in
   the former case.  Here any deviation from the
     true boundary implies a suboptimal exercise policy.
 
   We hasten to add 
    that the algorithm presented so far has one shortcoming
      related to the combined effects of finite sample size $N$ and  
      well-determined  values of the underlying variables at
     start time. 
    Only with  paths crossing the  YES-NO boundary its optimal location
      can be determined, of course. 
    Formally, the proof of proposition P1 requires that the probability 
      densities appearing 
    in eq.~(\ref{eq_optim_cont}) be positive-definite in a region
     covering  the YES-NO boundary.
       However,  this condition will be violated
     in practice  at early times, since all sampled paths
     emanate from a common starting point.
    On the other side, in such a situation the knowledge of the precise location 
     of the boundary becomes irrelevant for the valuation problem. It suffices
    to know at this stage on which side of the  boundary the vast
     majority of  paths can be found. 
    Two variants with differing degree of sophistication 
     have been developed  to determine the earliest time
    for which the algorithm above is to be used:
  \begin{description}
   \item[(7a)] 
     If  the  YES-NO boundary  gets closest to
       either one of  the endpoints
       of  $H$, e.g.\ 
       the paths with smallest or largest $S_t$ value 
        for plain vanilla options, 
      it is assumed that the boundary can be found outside the
       range of the path sample for all earlier times.
   \item[(7b)] 
    ``Flashlight mechanism'': 
      when the  YES-NO boundary moves towards a region only poorly
     covered by the original path sample, additional path segments are created
      which evolve in the region of the  YES-NO boundary at this time. 
    As far as  the proof of  proposition P1 is concerned
    the  probability   densities  in eq.~(\ref{eq_optim_cont})
      do {\em not } need to be consistent with the initial value
      specification for the pricing problem.    
    Here  we  make use of this freedom to generate a set
      of paths whose initial values are tuned to `shed light' on the
     location of the  YES-NO boundary. 
    
  \end{description} 
   Concerning accuracy
   simulation mode {\bf (7b) } is only barely noticably superior over mode
     {\bf (7a) } for the tested sample of plain vanilla options. 
   All simulation results 
    for more complicated options reported on in this paper 
     have been achieved employing mode {\bf (7b) }.
        
\begin{figure}[h]

\centerline{\hbox{
\psfig{figure=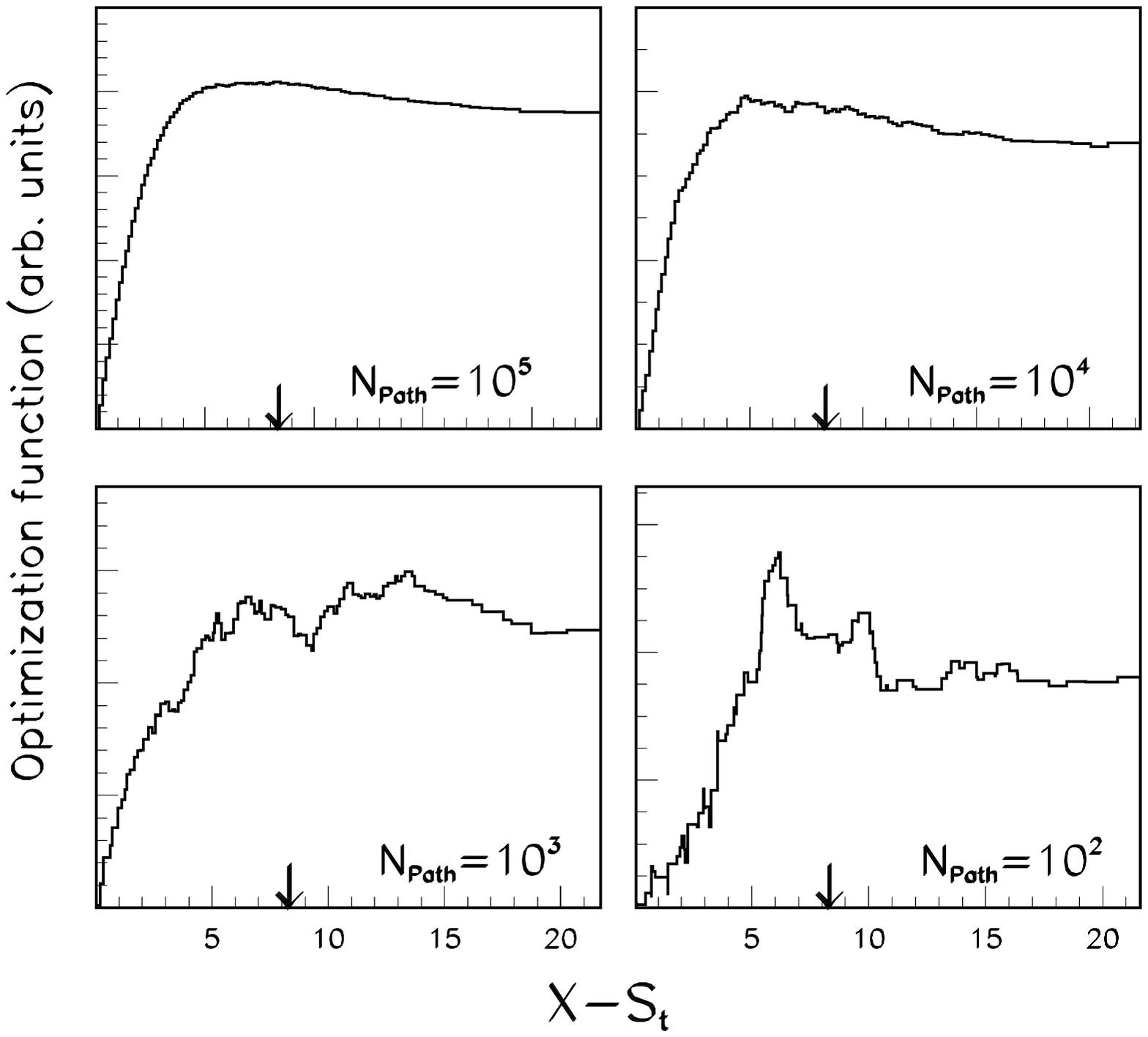,width=15cm,height=15cm}}}

\caption
[
 ]
{
 \label{fig_1}
}
\end{figure}

   At finite $N$ 
   the values and therefore also the  maximum of  the payoff function
   ${\cal O} _t^N ( S') $  
    do not change 
    between  values  of its argument 
    which are taken by two neighboring paths.
   Furthermore, 
    the  precision with which the true maximum 
    of function  ${ O} _t^\infty ( S') $
    can be estimated from a  finite size  sample 
     is  a  multiple ($\gg 1$) of the typical distance 
     between paths in the region along the   YES-NO boundary.
   Fig.~\ref{fig_1}  shows
 the Monte Carlo estimation of the optimization function
  $ {\cal O} _t^N ( S')$ for a plain vanilla put option
  on a security of (time-dependent)
  value  $S$
  in order to illustrate the finite $N$ effect.
  Selected  $ {\cal O} _t^N ( S')$  values
   are displayed at a fixed  time but varying the size of the
   path sample.
   The true minimum is quite shallow  for the selected time
   as can be  seen from the figure.
   The shallowness is related to the strong time dependence 
    of the location of the YES-NO boundary, because we consider 
   a time shortly before expiry. 
  The location of the maximum can be reliably extracted from the
   simulations only for
   rather large  samples on the order of $10^5$ paths.
  One may ask whether in view of the imprecisions 
   associated with any  finite sample size it is worthwile to search for the 
   exact maximum of  $ {\cal O} _t^N ( S')$. In fact,
    a well-chosen `smoothing' procedure which subtracts the white random
   noise may give more accurate results than taking the 
    real maximum of   $ {\cal O} _t^N ( S')$ as a boundary point.
   For this paper we compared  two variants how to determine 
    the  YES-NO boundary points:
  \begin{description}
   \item[(3a)] 
   Determine the  YES-NO boundary point as the element of
   the set of all intersection points $S_i(p)$ of the
   sampled paths with path $H$ for which ${\cal O} _t^N $
      restricted to these points 
     attains  a maximum value  $\max _{ \{S_i(p)\} } {\cal O} _t^N $.
   \item[(3b)] 
    Determine the  YES-NO boundary point as  the element of 
     the  set of all grid points $ S_i ^G $ for which ${\cal O} _t^N $  
      restricted to these points 
      attains a maximum value  $\max _{ \{ S_i(G)\} } {\cal O} _t^N $. 
    The distance between sites on the grid is lowered stepwise until
     some pre-defined precision is reached or
     more than one local extremum appears for the set of  values of 
       ${\cal O} _t^N $   restricted to the grid sites.
 \end{description} 
 Ultimately, in the limit of infinitely many paths
 both variants  may  give  the same correct result for
 the optimization of  ${\cal O} _t^N $.
  However, computation employing variant {\bf  (3b) } 
  is much quicker. It requires a number of operations growing
  only  linearly with $N$ while variant  {\bf  (3a) }  
   needs at least on the order of $N \log N$ operations.
 
\begin{figure}[h]

\centerline{\hbox{
\psfig{figure=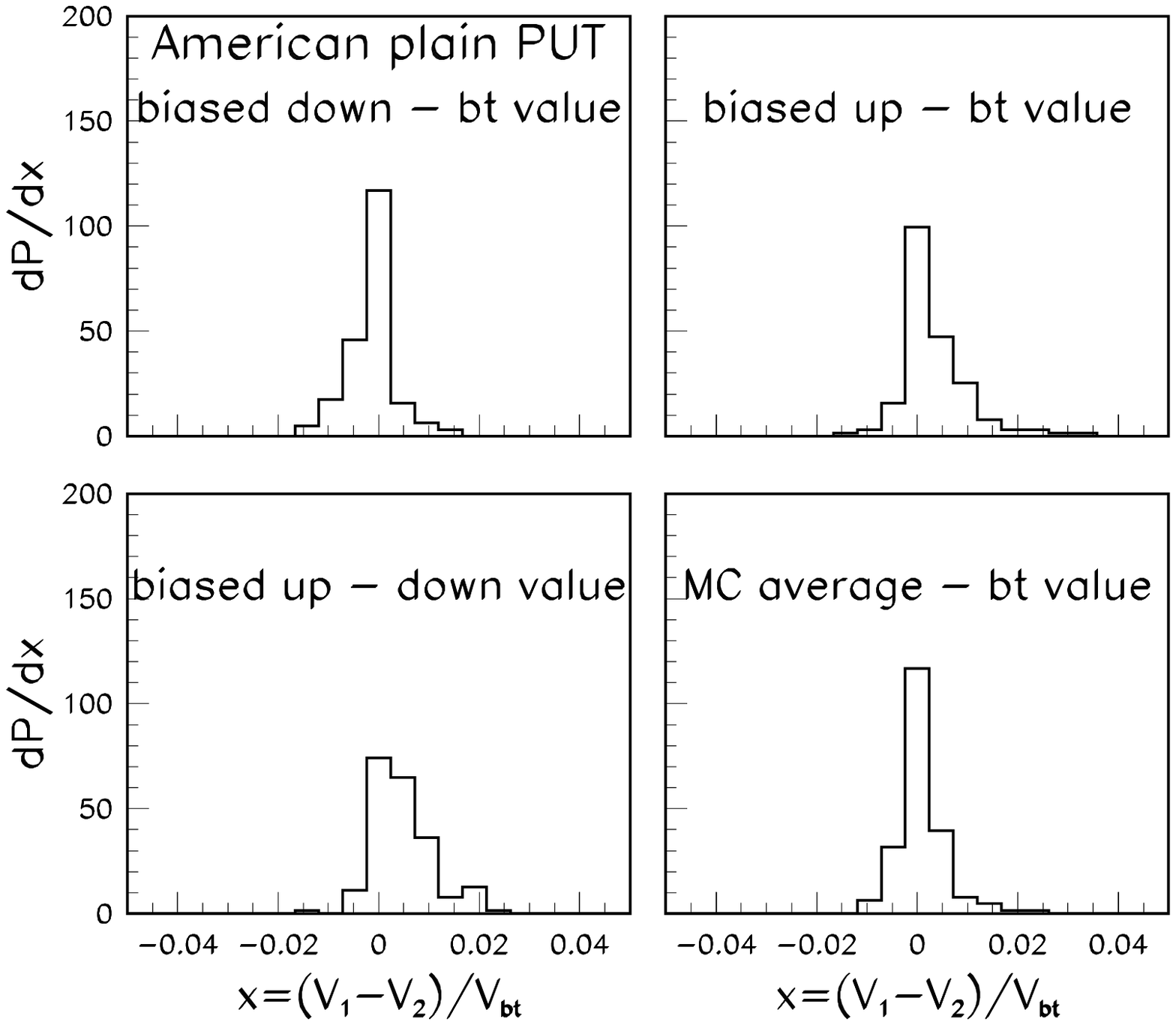,width=15cm,height=15cm}}}

\caption
[
 ]
{
 \label{fig_2}
}
\end{figure}

\begin{figure}[h]

\centerline{\hbox{
\psfig{figure=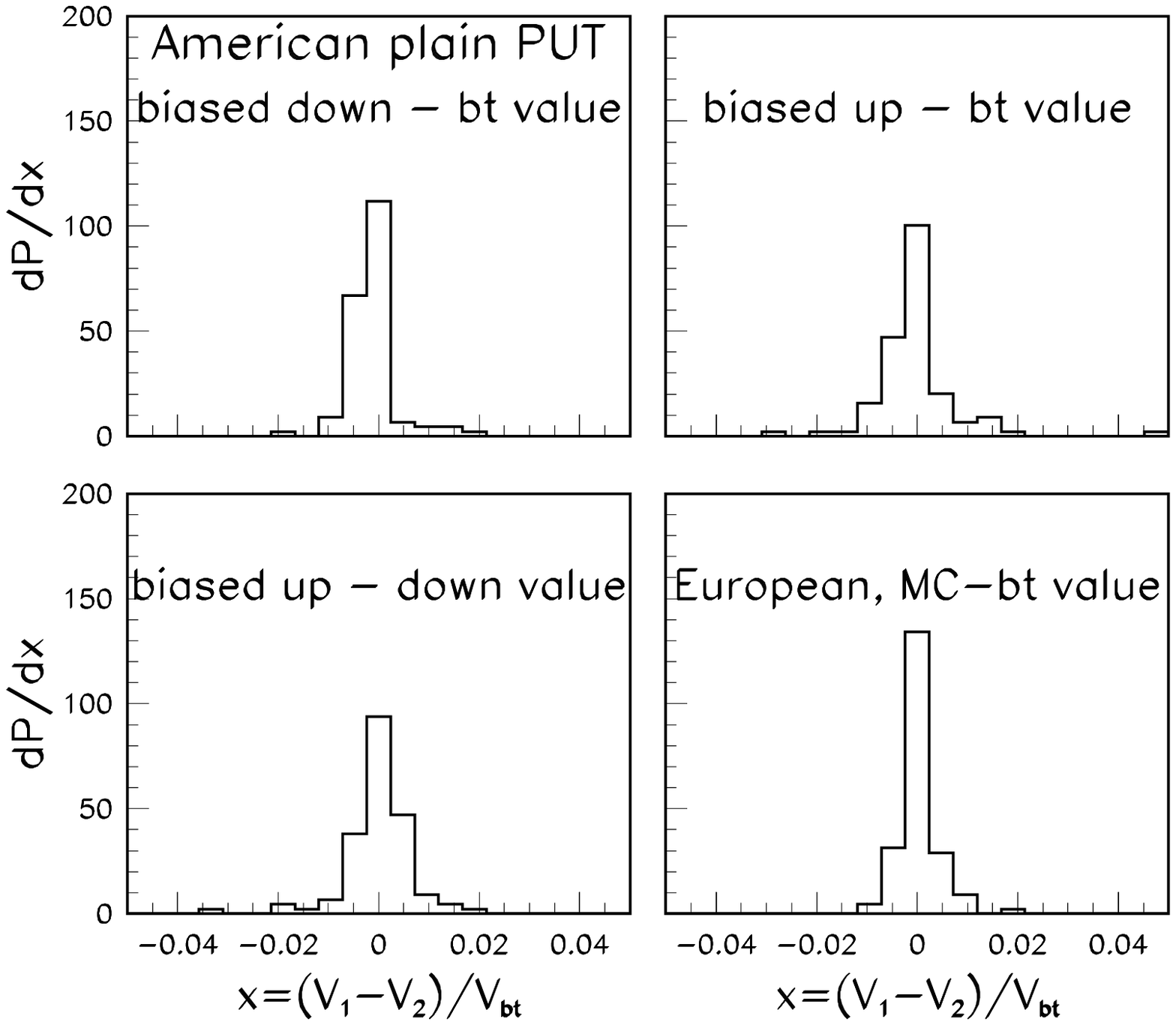,width=15cm,height=15cm}}}

\caption
[
 ]
{
 \label{fig_3}
}
\end{figure}

\newcounter{cnt_tab}
\setcounter{cnt_tab}{1}  
  We have generated Monte Carlo estimates for prices of  randomly sampled      
  standard American-style  put options 
  (see Table \Roman{cnt_tab} for details).
  We have employed the two simulation modes  {\bf (3a)  }  and  {\bf (3b)  } 
  to compare their effectiveness.
 The frequency distribution $dP/dx$ of
  the relative errors $x$ in   the simulation results based on
   Monte Carlo mode  {\bf (3a)  } 
    are shown in  Fig.~\ref{fig_2},
     for mode {\bf (3b)  } in  Fig.~\ref{fig_3}.
  These errors are determined from the normalized difference
   to the prices 
    extracted from  binomial-tree calculations:
 \begin{displaymath}
     x = ( V_0^{MC} -  V_0^{bt} )/  V_0^{bt}  \quad .
 \end{displaymath}
  The  Monte Carlo estimates have been calculated for
   the path sample which has been used to determine the
  YES-NO boundary and for an independently
  generated path sample. 
  The comparison between the results in the two modes reveals
  that both lead to similar accuracy. 
  In addition, the figures display the distribution of
   the differences between 
  Monte Carlo values  based on the original and the independent 
   path sample. 
   Upward and downward bias are clearly revealed for the
   simulation results in mode   {\bf (3a)  }
    (lower left panel of  Fig.~\ref{fig_2}).
  On the other side, the fact that 
    mode {\bf (3b)  } employs only an approximate optimization of
   the payoff function   ${\cal O} _t^N $  tends to wash out 
   most of the up- and downward biases.
  Fig.~\ref{fig_3} shows also the  error  distribution of  the corresponding
   European option estimations. 
  We note that the typical size of the errors is  comparable
   to the results for American puts. The error in both cases is thus
  dominated by the finite size of the path sample 
  ($N=10^5$ in  Figs.~\ref{fig_2} and \ref{fig_3}).
  The error induced by the bias 
   in the construction of the YES-NO  boundary at finite $N$ 
   seems to be of less importance.
  
  Let us turn now to the case of put options contingent on
  the geometric or arithmetic mean value of the underlying variable.
  As discussed in the preceeding section 
   exercise and hold regions and therefore the YES-NO boundary 
   are simply connected in the $S_t-\bar{S_t}(p_<)$ space. 
   The boundary is therefore one-dimensional in this space. 
   We have chosen a very simple procedure in order to track the
    boundary. These options are not directly contingent on the
    current $S$ value at any of the decision times which makes
    $\bar{S_t}(p_<)$ the more relevant variable.
   Therefore we  sort the sampled paths into bins covering the
     relevant $S_t$ region. Ignoring the differences between
    $S_t$ values {\em inside } each bin we 
     have reduced the optimization problem again to 
      finding the  $\bar{S_t}(p_<)$  value at which the
    payoff function  ${\cal O} _t^N $  attains maximum value. 
     This optimization proceeds as described already for plain
     vanilla options. 
    For the calculations described here we have taken the number
     of bins in  $S_t$  direction to be 20.
   
\begin{figure}[h]

\centerline{\hbox{
\psfig{figure=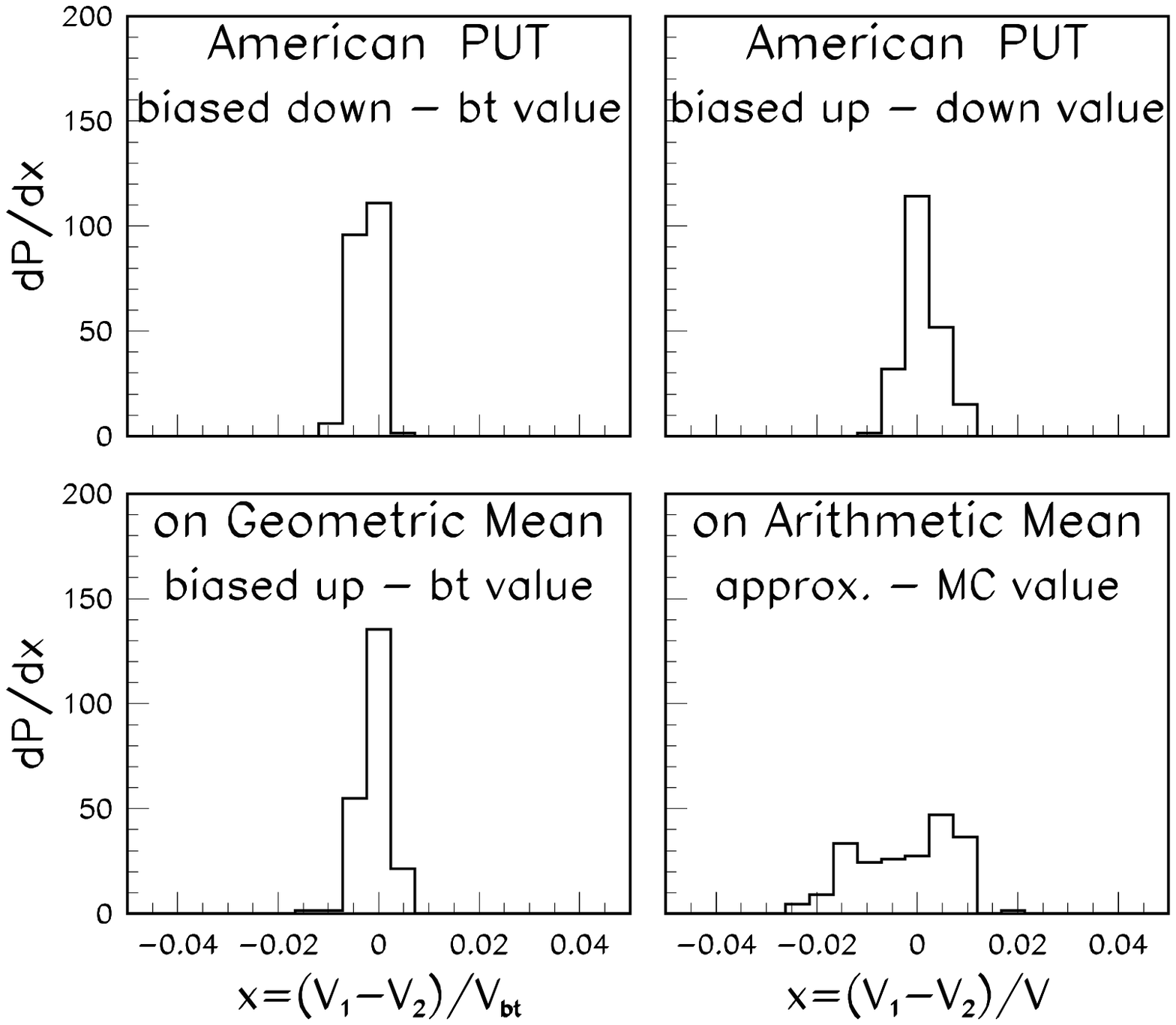,width=15cm,height=15cm}}}

\caption
[
 ]
{
 \label{fig_4}
}
\end{figure}

  Monte Carlo simulation
  results  for a randomly generated sample of put options 
  contingent on geometric averages 
     are presented in Fig.~\ref{fig_4}.
  As before, errors in the simulation have been calculated 
   in case of geometric averaging by comparing
  the Monte Carlo estimations 
   for the option prices to binomial tree results. 
  Here we consider only averages over the whole lifetime of the option.
  This is an $n^3$ problem for the binomial-tree method, $n$ being
   the number of time steps. Allowing the time window for the averaging
  to  slide would increase the storage and CPU time requirements
     by one power of $n$.
  In contrast, the Monte  Carlo  approach is not impacted severely by
  such a generalization.
   
  A comparison with the corresponding errors of  the
 simulation results for European options under the same conditions
   (not shown)
  reveals again that the errors are dominated by the finite size
  of the path sample and not by the uncertainty in determining the
   YES-NO boundary.
  The effect of upward and downward bias in the estimations due to 
  choosing either  the  original or the independently generated path sample
   are very weak as can be seen from the lower left panel of
  Fig.~\ref{fig_4}.  
  As  explained before this is a direct result of the approximations
  in the  optimization procedure {\bf (3b)  }.
  It is noteworthy that the 
    bias from parametrizing the  YES-NO boundary on a grid 
       (in $S_t$ direction) is stronger than either
   bias from taking  one of the 
   two path samples. 
   The representation of  the  YES-NO boundary on a grid tends to 
   {\em underestimate } the option prices because of the coarse
   graining procedure involved.
   The small nonstatistical downward net bias can be read off from
   the average  deviation between
   Monte Carlo estimation and  binomial tree results for the
  options contingent on geometric averages.
  It amounts to -~0.24 \% for the results based on the original 
   path sample and  -~0.1  \%  for the independent  path sample. 
     
  Results for prices of options contingent on arithmetic averages
  over the whole lifetime of the option are also included in
  Fig.~\ref{fig_4} (lower right panel).  
  Since no other  accurate method is available, we present the
  distribution of 
  relative deviations between  the 
   values based on an approximate formula suggested
  by  Ritchken, Sankarasubramanian,  and Vijh (1989)
  and the simulation results.
  The approximate formula relates the  price  for a claim contingent
  on  the   arithmetic   average to the corresponding option price
   using geometric averaging via 
 \begin{equation}
  \label{eq_approx}
   A_{am} \approx
         A_{gm}^{bt} - E_{gm}^{bt} + E_{am}^{MC}
        \quad  .
 \end{equation}    
  Here $A$ ($E$) denotes ``value of American (European) option''. 
  The subscript $am$ ($gm$) refers to 
  arithmetic (geometric) averaging.
  The superscripts $bt$ and $MC$ characterize the method of 
  calculation, binomial tree and Monte Carlo respectively.
   Fig.~\ref{fig_4} shows that the approximation
  works reasonably well.

\setcounter{cnt_sect}{3}
\section*{ \Roman{cnt_sect}  Conclusions}
    
 We have shown in this paper how to construct Monte Carlo algorithms for 
  American option valuation. 
   We are in the lucky situation that the suggested 
   algorithm turns out to be conceptually  simpler
    than the other   algorithms 
    suggested so far. 
    It depends only linearly on the number of sampled paths and 
    exercise dates which makes it  computationally feasable to 
    get close to the continuum limit   for these two variables. 
    Thus no price is to be  paid for the  accuracy 
    of this novel  Monte Carlo algorithm 
    in terms of  more CPU time or 
     exceedingly large storage requirements. 
    The crucial element, the optimization of 
        a certain  payoff function of the path sample, turns out to be
      numerically rather stable. 
    The tracking of the YES-NO boundary is achieved most precisely 
       in those regions of the state space covered  densely by
    sampled paths, i.e.\ where it is needed
     for the option price estimation. There is a further reason for the
     algorithm's stability. In case 
     a maximum turns out to be shallow 
      and is therefore not easily found by the algorithm, the option
      price does not depend much on the
     exact location of the  boundary between early-exercise and hold region.
    Indeed, we have demonstrated that the errors in the simulation results
    are typically dominated by  the statistical errors if  the  size of
    the path sample is on the order of $10^5$. 
    This points to one important direction of  future research.
    Paskov and Traub (1995)   have  shown  
    that the replacement of pseudo-random numbers
    in the traditional Monte Carlo approach by a series of so-called
    quasi-random numbers may lead to a spectacular gain  in the
   precision of the achievable results, at least
   for some problems. The reason is
   that the   quasi-random numbers are distributed more evenly.
   By using appropriately chosen  sequences of quasi-random numbers
   sample sizes  for estimation of  European  options 
   could be reduced  by orders of magnitude, without loss of accuracy.
   It would open new dimensions for the applications of the novel
   Monte Carlo   algorithm  if this would also hold true for
   American  options.
  Moreover, we would like to add that the extraction of sensitivity
   coefficients -- usually called the ``Greeks''  -- 
   employing the novel Monte Carlo algorithm is as straightforward as
   for  previously considered situations, e.g.\ by 
  Fu and Hu  (1995). 
      
  We have presented some simulations for option types 
  for which other means of calculation are known.       
  However, the  Monte Carlo   algorithm presented here provides
   unique opportunities to evaluate path-dependent option prices
   for which no other methods are available or computationally feasible.
   One example, 
  the case of options contingent on arithmetic averages has been discussed
  already  in this paper. 
   Another  candidate are ``look-back'' options of American style.
  Usually, look-back   options allow the owner to buy or sell an asset
   for a price dependent on the values of the asset during the whole
  lifetime until expiry,  in the simplest case the minimum or maximum.
  Restricting the time window for the look back
   but allowing for early exercise  adds  American features to these
   rather common options.
  Another virtue of the  presented  Monte Carlo   algorithm is to
   allow  option  valuation for more complicated stochastic processes
  than Wiener-type diffusion. 
      
\newpage     
      
 {\large \bf References: }
 \vspace{1.0cm} 

{\noindent 
 J.~Barraquand and D.~Martineau, 1995, 
 Numerical valuation of High Dimensional  Multivariate 
  American Securities,
  {\em Journal of Financial and Quantitative Analysis} 30, 383-405.

 F.~Black, and M.~Scholes, 1973,
 The Pricing of  Options  and Corporate Liabilities,
  {\em Journal of Political Economy} 81, 637-654.

 P.~Boyle, 1977, 
  Options: A Monte Carlo Approach,
    {\em Journal of Financial Economics} 4 , 323-338.

 P.~Boyle, M.~Broadie, and P.~Glassermann, 1997,
  Monte Carlo Methods for Security Pricing,
  in  Risk Publication Group, ed.: 
   Pre-Course Material 
  for RISK$^{TM}$ Training Course (January 1997)
   `Monte Carlo Techniques For effective Risk Management And
   Option Pricing', (Risk Publication Group, London).
    
 M.~Broadie, and P.~Glassermann,  1995,
  Pricing American-Style Securities Using Simulation,
  Working paper  (Columbia University), to appear in
   {\em Journal of  Economic Dynamics and Control }.

 J.C.~Cox, S.A.~Ross, and M.~Rubinstein, 1979, 
 Option Pricing: A Simplified Approach, 
  {\em Journal of Financial Economics} 7, 229-263.

M.C.~Fu and J.Q.~Hu, 1995,  
 Sensitivity Analysis For Monte Carlo Simulation  Of
  Option Pricing,  
   {\em Probability in the Engineering and
 Informational Sciences} 9, 417-446.

  R.~Geske, 1979, 
  A Note on an  Analytic Valuation Formula for  Unprotected  
  American  Call  Options  with Known Dividends, 
    {\em Journal of Financial Economics} 7, 375-380.

  D.~Heath, R.~Jarrow and A.~Morton, 1992, 
  Bond Pricing and the Term Structure of Interest Rates:
    A New Methodology,
    {\em Econometrica } 60, 77-105.

 J.C.~Hull, 1993.
 Options, futures, and other derivative securities,
  (Prentice-Hall, Englewood Cliffs).

 R.C.~Merton, 1973, 
 Theory of Rational  Option  Pricing,
   {\em Bell  Journal of Economics and  Management Science}
     4, 141-183.

 S.~Paskov, and J.~Traub, 1995,
 Faster Valuation of financial Derivatives,
    {\em Journal of Portfolio Management }  22, 113-120.

P.~Ritchken, L.~Sankarasubramanian, and A.M.~Vijh,  1989, 
 The Valuation of Path Dependent Contracts on the
  Average, Working Paper, Case Western Reserve University
  and  University of Southern California.

  R.~Roll,  1977, 
  An Analytic Formula for Unprotected  American  Call  Options
    on Stocks with Known Dividends,
    {\em Journal of Financial Economics} 5, 251-258.

 J.A.~Tilley, 1993, 
 Valuing  American Options  in a Path Simulation Model,
  Transactions of the Society of Actuaries 45, 83-104.
    
  R.~Whaley, 1981, 
  On the  Valuation of  American  Call  Options 
     on Stocks with Known Dividends,    
     {\em Journal of Financial Economics} 9, 207-211.

}

\newpage
 
{\noindent \bf { \large Table \Roman{cnt_tab}: } 
\begin{minipage}[t]{12.0cm}
 Parameters for random  ideal \\ 
  Wiener processes and   put options 
\end{minipage}         
               }
\vspace{2em}     

 \begin{tabular}{|c||r|r|r|r|r|}
                                  \hline  
   Parameters: &  $r$   &  $\sigma $   &  $S_0$  &  X  &   T   \\
                                \hline      \hline  
   mean value &  0.10 & 0.40 & 100.0  & 100.0 & 0.50 \\
                                  \hline  
   standard deviation &  0.05 & 0.20 &  50.0  &  50.0 & 0.25   \\ 
                                  \hline  
 \end{tabular} 
\vspace{2.2cm}

{\noindent        
   Caption:
   Mean values  of riskless interest rate $r$, volatility  $\sigma $,
   start price  $S_0$,  strike price $X$ and time until expiry $T$
   are given in first row, their
   standard deviations respectively in second row.
   The random values are chosen according to 
   independent  normal distributions.
   In a second step, 
   some parameter values may be rejected and replaced by other random 
   values if  they are nonpositive. 
   Furthermore,  strike price $X$ is required to be within 
   a 2 (1)  $\cdot \sigma \sqrt{T} $  range around  the 
   mean $S_T$ value at maturity  $S_0 + r \cdot T$ for
   options contingent  on the current  (average) value. 
   This constraint could be relaxed using a more sophisticated
   importance sampling scheme.
   The number of evenly spaced time steps $N_t$ for the
   Monte Carlo simulation and the binomial-tree calculations
   has been fixed arbitrarily to be equal 100.
}
      
\newpage
 
{\noindent  \LARGE   Figure Captions:}
\vspace{0.6cm}

{\noindent \large Figure 1: }

{\noindent        
 Monte Carlo estimation of the optimization function
  $ {\cal O} _t^N ( S')$ for a plain vanilla put option.
  Option and stochastic process parameters are chosen
  as in Hull's book (1993), Example 14.1. 
  The Monte Carlo simulation with number of time steps equal 50
  has been repeated  varying the total number of sampled
  paths: $10^2$, $10^3$, $10^4$ and $10^5$.
  Only the function values 
  at   time step 45  are displayed  
  for 100  values of the argument  $S'$.
  The  location of the YES-NO boundary according to a
  binomial-tree calculation is pointed at by an arrow.
  
}
\vspace{0.2cm}

{\noindent \large Figure 2: }

{\noindent        
 Frequency distribution of
  the relative errors in 
  simulation estimates of a sample of plain put option prices.
  The error is determined from the normalized difference
   to the prices $V_{bt}$
    extracted from  binomial-tree calculations.
  The  Monte Carlo estimates are calculated for
   the path sample which has been used to determine the
  YES-NO boundary (biased up) as well as for an independently
  generated path sample (biased down) and the average of the 
   two samples (MC average).
  In addition, the figure displays the distribution of
   the difference between upward and downward biased 
  Monte Carlo values (lower left  panel).
  Each path sample encompasses  $10^5$ generated paths
   per option.
  The  Monte Carlo simulation employed the mode
  in which the `global maximum' of the optimization function
   $ {\cal O} _t^N ( S')$ was searched. 
  
}
\vspace{0.2cm}

{\noindent \large Figure 3: }

{\noindent        
 The content of this figure is the same as in the previous figure,
 with the exception of two features.
 The mode of the  Monte Carlo simulation has been to look
  up the maximum of  the optimization function 
   on a grid (see text). 
   Furthermore, the distribution of
 errors for the averaged Monte Carlo values 
   which looks very similar to the corresponding distributions
   of upward and downward biased estimates
  is not displayed. Instead, 
   the distribution of the errors for the corresponding
   {\em European } option values is displayed 
  (in the lower right panel).
}
\vspace{0.2cm}

\newpage

{\noindent \large Figure 4: }

{\noindent        
 Frequency distribution of
  deviations between
  simulation estimates of  put option prices
   on geometric (left side)
   and arithmetic average values (right side)
   and benchmark calculations.
  The benchmark for the prices on geometric
   averages is provided by the 
   rather accurate  binomial-tree  method, 
   and for the  arithmetic averages by
    the approximate formula eq.~(\ref{eq_approx}). 
}
\vspace{0.2cm}     

\newpage

\newpage

\newpage

\newpage

\newpage

\end{document}